# Preliminary Functional-Structural Modeling on Poplar (*Salicaceae*)

Dongxiang Liu[1], Mengzhen Kang[2,*], Véronique Letort[3], Meijun Xing[1], Gang Yang[1], Xinyuan Huang[1], Weiqun Cao[1]

[1]*School of Information Science & Technology, Beijing Forestry University, Beijing, China*
[2]*LIAMA, Institute of Automation Chinese Academy of Science, Beijing, China*
[3]*Laboratory of Applied Mathematic Ecole Centrale Paris, Paris, France*
[*]Corresponding author *(mzkang@nlpr.ia.ac.cn)*

**Abstract**

*Poplar is one of the best fast-growing trees in the world, widely used for windbreak and wood product. Although architecture of poplar has direct impact on its applications, it has not been descried in previous poplar models, probably because of the difficulties raised by measurement, data processing and parameterization. In this paper, the functional-structural model GreenLab is calibrated by using poplar data of 3, 4, 5, 6 years old. The data was acquired by simplifying measurement. The architecture was also simplified by classifying the branches into several types (physiological age) using clustering analysis, which decrease the number of parameters. By multi-fitting the sampled data of each tree, the model parameters were identified and the plant architectures at different tree ages were simulated.*

## 1. Introduction

Poplar is a member of *Salicaceae* family and is one of the main species for windbreaks and timber in plain and sand [1]. It includes several species and is distributed in the world extensively. It is cultivated in plantations not only for paper, wood products and energy, but also for soil erosion control [1]. In addition, it has proven to be efficient for alleviating environmental toxins pollution [2].

As poplar's wide range of use, many studies on modeling its growth have already been carried out. Such as plant physiology, morphology and genetics. Hao *el al.* [3] built a model for growth and harvest prediction for the poplar forest; Rauscher *el al.* simulate the ecophysiological growth processes of poplar seedlings with the ECOPHYS [4] [5]; and zhang *el al.* use the Penman-Monteith model to simulate the transpiration volume of poplar at growing season [6].

These models are helpful to understand crop production and management, but as poplar's detailed architecture is not described, it is difficult to take into account the influence of the plant structure on wood quantity and quality.

The plant architecture is increasingly recognized to be a determinant factor of plant growth and it is necessary to take it into account in plant growth modeling [7]. This boosts the study of functional-structural models [8] [9]. Functional-structural models are potentially able to simulate the interactions between the architectural development of plants and their physiological functions. Recently, several functional-structural plant models have been put forward, such as L-Peach [10] and LIGNUM [11]. These models can simulate the growth of plants, but the computing complexity, computing accuracy and parameter estimation become a huge obstacle when using these models to simulate complex trees with a large number of organs [12].

In the past years, a functional-structural plant model GreenLab was developed. Its parameters are estimated through the generalized least squares method [13]. It allows fitting all the model parameters according to target data of different kinds of organs collected on plants observed at different growth stages. This is called multi-fitting technique [14]. Multi-fitting technique can minimize the error between model output and the target data of different stages [14]. When modeling trees with GreenLab, the topology of trees can be simplified based on the notion of substructure, which makes it easy to deal with trees with complex architecture. Applications of GreenLab have been made on Chinese pine sapling [15] and beech trees [12].

In this paper, we apply GreenLab to simulate the growth of poplar through model calibration on data of 3, 4, 5, 6 year-old trees. We introduce briefly the calibration process including data acquisition, data simplification, fitting target and simulation result.

## 2. Material and method

### 2.1. Plant material and sampling

The Poplar variety measured is '107'. Samplings were conducted in June 2008 at Fugou country, Henan province (34°3′N, 114°23′E). Four trees that had grown 1, 2, 3, 4 years were destructively measured. They are noted tree1, tree2, tree3 and tree4. These trees were generally plugged in spring and implanted to the nursery in autumn. The next March, they were replanted to the field where the sampling was made, so their ages are actually three, four, five, six years old. The planting density in field is $3m \times 4m$.

There are a lot of branches with similar structure inside a growth unit (GU) of these trees. As measuring all the branches in details would be a very tedious work, the measuring was simplified by sampling one or more representative branches according to their length and morphological characteristics. In our experiment, tree1 and tree3 are respectively collected four branches samples, the tree2 and tree4 are collected six.

For each sample branch axis, the fresh weight, length and diameter of each internode and the number of leaf scars of each GU were recorded. Three leaves were sampled at different position of each GU, and the fresh weight and surface were recorded. The fresh weight was weighted by a balance (accuracy is 0.001g) and an electronic scale (accuracy is 1g). The length was measured with a ruler (accuracy is 0.01cm) or tape (accuracy is 0.1cm). The diameter was measured with a vernier caliper (accuracy is 0.02mm). The leaf surface area was obtained by scanning the leaf and using image processing software to count the pixels.

### 2.2 GreenLab model

GreenLab model has been presented in previous papers [8] [16], and a tutorial is available on the website (www.greenscilab.org). In our study, only the deterministic version of the model was used. That is to say we don't consider the stochastic and environment variations. The multi-fitting method [14] was used to optimize parameters in our study, hence permitting simultaneously optimizing all relevant parameters against many target files of botanical and morphological plant observations, and being able to construct identical phenotypes by implementing the same rules and parameters.

## 3. Results

### 3.1 Classification of branch physiological age

In the target file, we define each axis in accordance with the physiological age (PA) and chronological age (CA), and describe it by GU with average internode weight, length and diameter and blade weigh and surface. The PA is the number of states for axis in a plant, which characterizes the potential evolution strength [12]. In this paper, the PA of the branch was determined by cluster analyzing of the weight of new internode at the top of the branch axis. And the result is accordance with the PA determined by characteristics and length of the branch [12].

### 3.2. Sink strength of internodes

In GreenLab hypothesis, there is a simple proportional relationship between the internode weight and the leaf weight of the new GU. We consider the sink of leaf in the main axis as a reference value, which was set to 1. So the internode sink of different PA is a relative value, and we can get the internode sink from the ratio between internode weight and leaf weight of the new GU. The *Figure1* shows this ratio. The result indicates that there is a strong linear relationship between leaf mass and internode mass of each PA. And the result is listed in *Table1*.

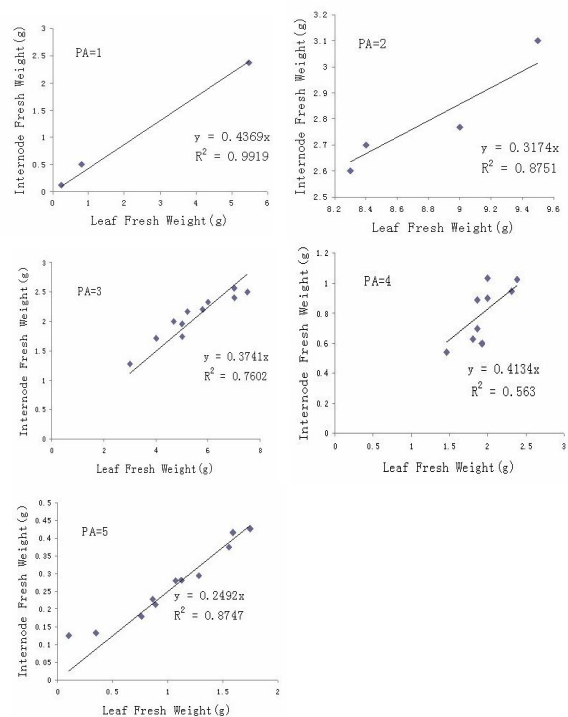

**Figure** 1**. The ratio between the leaf fresh weight and the internode fresh weight**

**Table 1. Sinks and allometric parameters**

| Parameters | | values | | | | |
|---|---|---|---|---|---|---|
| | | PA1 | PA 2 | PA3 | PA4 | PA5 |
| $p$ | Internode sink | 0.561 | 0.726 | 0.858 | 0.954 | 0.499 |
| $b$ | Int allometry | 5.193 | 7.229 | 9.279 | 12.146 | 20.78 |
| $\beta$ | Int allometry | -0.22 | -0.5512 | -0.18 | 0.0576 | 0.259 |
| $\varepsilon$ | SLW | 0.0287(g/cm2) | | | | |

### 3.3. Internode allometry and specific leaf weight

In GreenLab model, the organ shape can be determined by using simple allometric rules based on organ weight, length and diameter. The parameters $b$ and $\beta$ determine the internode shape. They can be estimated by analyzing the weight and length of new internode located at top of branch axis of each PA. We can see from *Figure 2* that the internode length and internode weight have a strong relationship. The estimated results of internode allometric parameters are listed in *Table1*.

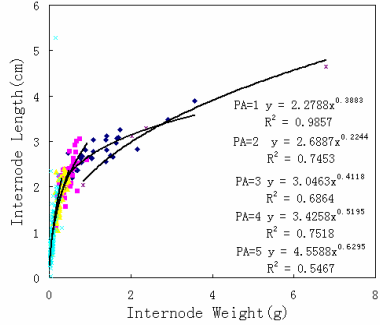

**Figure 2. Internode allometric parameters**

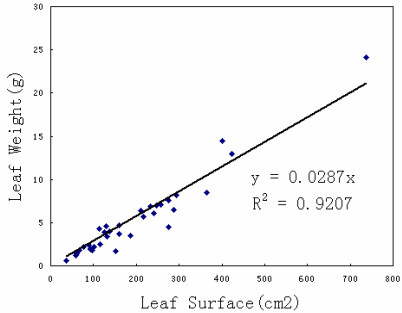

**Figure 3. The specific leaf weight**

The specific leaf weight (SLW) is varied from trees of different age, but the data of each tree we sampled are so few that we can't get a good value of the SLW of each tree. So we set the SLW of each tree as a same constant in our study, which we get from the ratio between the leaf weight and the leaf surface of all blade data in *Figure 3*.

### 3.4. Fitting the target data

The hidden parameters of GreenLab are estimated by calibration of the model on target data. In this paper, the data of four trees were multi-fitted to enhance the accuracy of parameter estimation. Trees share the same functional parameters (except seed biomass) but their topology is fixed as recorded in the target files.

**Table 2. Hidden parameters**

| Parameters | | Values | | | |
|---|---|---|---|---|---|
| | | tree1 | tree2 | tree3 | tree4 |
| $Q_0$ | seed biomass | 14.23 | 0.97 | 47.23 | 10.94 |
| $r_p$ | resistance | 6.4319 | | | |
| $p_c$ | ring sink | 0.13882 | | | |

Based the measured data of poplar, some parameters can be directly estimated, so only 3 parameters need to be fitted in our study: the initial seed biomass $Q_o$, the coefficient for leaf resistance $r_p$ and the relative sink strength of ring demand $p_c$. Calibration results are shown in *Table 2*. *Figure 4* illustrates the comparison between simulated and measured data for some of the outputs, such as the internode diameter of each GU (4a), the internode mass of each GU (4b), cumulated internode mass (4c) and cumulated blade mass (4d). Although the fitting is imperfect, the global trend is correct.

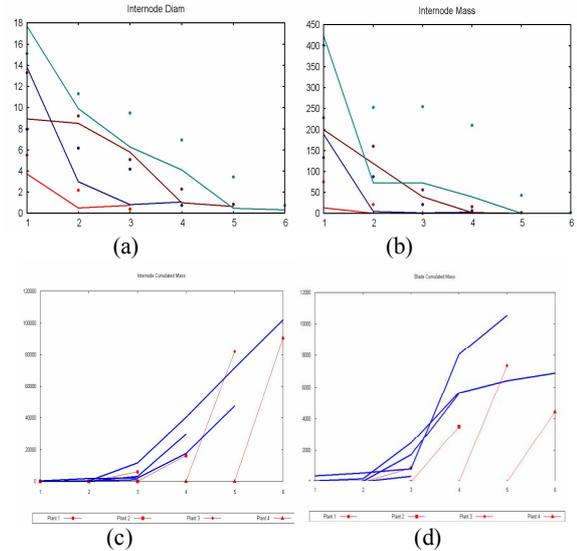

**Figure 4. Multi-fitting results**

### 3.5. Topological structure of plants

The topological structures of the four trees are shown in *Figure 5*. The positions of branches were not recorded and are arbitrarily set but the number of organs corresponds to the observations.

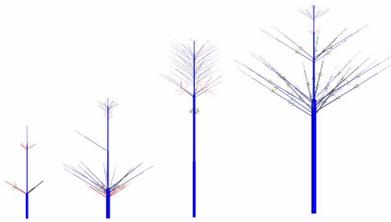

**Figure 5. The topological structures of four trees, of age 3, 4, 5, 6 years old.**

## 4. Discussion and Conclusion

The functional-structural model GreenLab has been widely used in crops, but only preliminary used in trees, such as pine and beech tree. This paper is an exploratory study on poplar using GreenLab. It analyses the ratio of leaf fresh weight to internode with experimental data, estimates the model parameters by model inversion, renders the simplified topology, and completes trees simulation.

Although the identified model parameters are approximately consistent with the results, further potential improvements can be achieved. First, we need to get more complete data to further work. Second, the sink dynamics of root system would be got better definition by collecting appropriate data. Third, the environmental effects should be considered in further studies as our long-term objective is the integration of environmental factors on the growth of poplar. Fourth, in the simulation, the distribution of branches should be set according to observed position of branch.

This work was the first experience of analyzing poplar growth with GreenLab. It has highlighted several weak points of the model but it also brings experience for the future works on more complex trees.

## Acknowledgments

We particularly thank the students and Henan Bureau of Forestry involved in collecting the experimental data. This study was funded by National Hi-Tech 863 Program of China (No.2006AA10Z232), Natural Science Foundation of China (No.60703006), the Excellent SKL Project of NSFC (No.60723005) and the National Science and Technology supported project (No.2006BAD10A03).